\documentclass[a4]{amsart}
\usepackage{amssymb,latexsym}
\begin{document}
\newtheorem{theorem}{Theorem}
\newtheorem{lemma}[theorem]{Lemma}
\newtheorem{sublemma}{Sublemma}
\newtheorem{proposition}[theorem]{Proposition}
\newtheorem{corollary}[theorem]{Corollary}
\renewcommand{\thefootnote}{\fnsymbol{footnote}}

\newcommand{\myfbox}{\fbox}  \renewcommand{\myfbox}{\relax}
\newcommand{\name}[3]{\addtocounter{num}{1}\hbox to 1.7cm{$#1#2_{#3}$\hfill} }
\newcounter{num}
\newcommand{\knum}[5]{\put(0,115){\name{#1}{#2}{#3}}}

\title[Prime knots with arc index 12 up to 16 crossings]%
{Prime knots with arc index 12 up to 16 crossings}

\author{Gyo Taek Jin}
\address{Department of Mathematical Sciences,
Korea Advanced Institute of Science and Technology,
Daejeon, 34141, Korea}
\email{trefoil@kaist.ac.kr}

\author{Hyuntae Kim}
\address{Department of Mathematical Sciences,
Korea Advanced Institute of Science and Technology,
Daejeon, 34141, Korea}

\author{Seungwoo Lee }
\address{Moasys Corporation, Ilkwang Bldg 4th, 220 Baumoe-ro Seocho-gu Seoul, 06746, Korea}

\author{Hun Joo Myung }
\address{Korea Institute of Science and Technology Information, 245 Daehak-ro Yuseong-gu Daejeon, 34141, Korea}

\begin{abstract}
As a continuation of the works \cite{Jin2006} and \cite{Jin2011}, we provide the list of prime knots with arc index 12 up to 16 crossings and their minimal grid diagrams.
There are 19,513 prime knots of arc index 12 up to 16 crossings.
\end{abstract}


\maketitle

\section{Introduction}
An \emph{arc presentation\/} of a knot is an embedding of a knot into the union of finitely many vertical half planes whose common boundary is the $z$-axis so that each of the half planes intersect the knot in a single properly embedded curve. These curves are called the \emph{arcs\/} of the arc presentation. The minimal number of arcs among all arc presentations of a given knot is called the \emph{arc index}~\cite{C1995}.

A \emph{grid diagram\/} is a knot diagram whose projection is a closed curve which is composed of finitely many horizontal line segments and the same number of vertical line segments such that vertical line segments always cross over horizontal line segments.
An $n\times n$  \emph{Cromwell matrix\/} is an $n\times n$ square matrix whose entries are either 0 or 1 such that there are exactly two 1's in each row and each column. By joining the 1's in a Cromwell matrix by horizontal line segments and vertical line segments, we obtain a grid diagram. See Figure~\ref{fig:grid&cromwell}.

On the other hand, a grid diagram can be converted easily to an arc presentation with the number of arcs equal to the number of vertical line segments as indicated in Figure~\ref{fig:grid2arc}.
Conversely, an arc presentation can be easily converted to a grid diagram with the number of vertical line segments equal to the number of arcs.

\begin{figure}[h]
\centering
\setlength{\unitlength}{0.6mm}
\begin{picture}(40,40)
\thicklines \put(0,0){\line(1,0){30}} \put(30,0){\line(0,1){30}}
\put(30,30){\line(-1,0){7}} \put(17,30){\line(-1,0){7}}
\put(10,30){\line(0,-1){20}} \put(10,10){\line(1,0){17}}
\put(33,10){\line(1,0){7}} \put(40,10){\line(0,1){30}}
\put(40,40){\line(-1,0){20}} \put(20,40){\line(0,-1){20}}
\put(20,20){\line(-1,0){7}} \put(7,20){\line(-1,0){7}}
\put(0,20){\line(0,-1){20}}
\end{picture}
\qquad
\normalsize
\begin{picture}(50,45)(-5,2)
\put(0,40){0}\put(10,40){0}\put(20,40){1}\put(30,40){0}\put(40,40){1}
\put(0,30){0}\put(10,30){1}\put(20,30){0}\put(30,30){1}\put(40,30){0}
\put(0,20){1}\put(10,20){0}\put(20,20){1}\put(30,20){0}\put(40,20){0}
\put(0,10){0}\put(10,10){1}\put(20,10){0}\put(30,10){0}\put(40,10){1}
\put(0,  0){1}\put(10,  0){0}\put(20,  0){0}\put(30,  0){1}\put(40,  0){0}
\put(-5,-2){\line(0,1){47}}\put(-5,-2){\line(1,0){2}}\put(-5,45){\line(1,0){2}}
\put(47,-2){\line(0,1){47}}\put(47,-2){\line(-1,0){2}}\put(47,45){\line(-1,0){2}}
\end{picture}
\caption{A grid diagram of a trefoil knot and the corresponding $5\times5$ Cromwell matrix}\label{fig:grid&cromwell}
\end{figure}

\begin{figure}[h]
\centering
\setlength{\unitlength}{0.7mm}
\begin{picture}(50,55)(-10,0)
{
\put(-10,10){\line(1,-1){10}}
\put(-10,10){\line(4,-1){7.5}} \put(2.5,6.875){\line(4,-1){27.5}}
\put(-10,20){\line(2,-1){7.5}} \put(2.5,13.75){\line(2,-1){7.5}}
\put(-10,20){\line(5,-1){7.5}} \qbezier(2.5,17.6)(5,17.2)(7.5,16.6)
                               \put(12.5,15.4){\line(5,-1){15}} \put(32.5,11.6){\line(5,-1){7.5}}
\put(-10,30){\line(1,-1){10}}
\put(-10,30){\line(3,-1){17.5}} \put(12.5,22.5){\line(3,-1){7.5}}
\put(-10,40){\line(2,-1){20}}
\put(-10,40){\line(4,-1){27.5}} \put(22.5,31.875){\line(4,-1){7.5}}
\put(-10,50){\line(5,-1){50}}
\put(-10,50){\line(3,-1){30}}
\put(30,0){\line(0,1){30}}
\put(10,30){\line(0,-1){20}}
\put(40,10){\line(0,1){30}}
\put(20,40){\line(0,-1){20}}
\put(0,20){\line(0,-1){20}}
}
\qbezier[50](-10,5)(-10,30)(-10,55)
\linethickness{0.3mm}
\put(0.7,0){\line(1,0){28.6}}
\put(29.3,30){\line(-1,0){6.3}}
\put(17,30){\line(-1,0){6.3}}
\put(10.7,10){\line(1,0){16.3}}
\put(33,10){\line(1,0){6.3}}
\put(39.3,40){\line(-1,0){18.6}}
\put(20.7,20){\line(-1,0){7.7}}
\put(7,20){\line(-1,0){6.3}}
\put(29.3,0){\line(0,1){30}}
\put(10.7,30){\line(0,-1){20}}
\put(39.3,10){\line(0,1){30}}
\put(20.7,40){\line(0,-1){20}}
\put(0.7,20){\line(0,-1){20}}
\end{picture}
\caption{Construction of an arc presentation from a grid diagram}\label{fig:grid2arc}
\end{figure}

For this work, we've generated  $12\times12$ Cromwell matrices so that they cover all prime knots with arc index 12. Using `knotscape' we identified the knots corresponding to these matrices. Then we removed those corresponding to the unknot, composite knots, links, and knots with arc index less than 12. The remaining knots have arc index 12. Finally we removed duplicates.

Table~1 shows the number of prime knots of given arc index and minimal crossing number. 
In the columns of the arc index 11 and 12, the  knots counted by the thirteen italicized numbers may have smaller crossing numbers but not smaller than 17.
The five numbers in boldface are newly determined by this work.

\begin{table}[h]
\small
\caption{Prime knots up to arc index 12 or  crossing number 13}
{
\newcommand{\vsp}{\vrule width0pt height 10pt depth.5pt} 
\begin{tabular}{|c|c|c|c|c|c|c|c|c|c|c|c|c|}
\hline
\vrule width0pt height 12pt depth22pt
\lower15pt\hbox{Crossings} \kern-25pt Arc index &
        \hbox to 7pt{\hfil 5\hfil}&
        \hbox to 7pt{\hfil 6\hfil}&
        \hbox to 7pt{\hfil 7\hfil}&
        \hbox to 7pt{\hfil 8\hfil}&
        \hbox to 7pt{\hfil 9\hfil}&
        \hbox to 7pt{\hfil {10}\hfil}&
        \hbox to 7pt{\hfil {11}\hfil}&
        \hbox to 7pt{\hfil {12}\hfil}&
        \hbox to 7pt{\hfil {13}\hfil}&
        \hbox to 7pt{\hfil {14}\hfil}&
        \hbox to 7pt{\hfil {15}\hfil}& Subtotal\\
\hline
\vsp  3&1 & & & & &  &  &  &  &  &  &  1\\
\hline
\vsp  4& &1 & & & &  &  &  &  &  &  &  1\\
\hline
\vsp  5& & &2 & & &  &  &  &  &  &  &  2\\
\hline
\vsp  6& & & &3 & &  &  &  &  &  &  &  3\\
\hline
\vsp  7& & & & &7 &  &  &  &  &  &  &  7\\
\hline
\vsp  8& & &1&2& &{18} &  &  &  &  &  &21\\
\hline
\vsp  9& & & &2&6&  &{41} &  &  &  &  & 49\\
\hline
\vsp 10& & & &1&9&32&  &  123 &  &  &  &  165\\
\hline
\vsp 11& & & & &4&46&135&  & 367 &  &  & 552\\
\hline
\vsp 12& & & & &2&48&  211&   627&  &1288 &  & 2176\\
\hline
\vsp 13& & & & & &49&  399 &\bf 1412 &\bf 3250 & &4878 &9988\\
\hline
\vsp 14& & & & & &17&  477  & \bf 3180\\
\cline{1-9}
\vsp 15& & & & &1&22&  441 & \bf 6216 \\
\cline{1-9}
\vsp 16& & & & & & 7&  345 & \bf 7955 \\
\cline{1-9}
\vsp 17& & & & & & 1&  192  &  10283\\
\cline{1-9}
\vsp 18& & & & & &  &  \it 75 & \it 8584\\
\cline{1-9}
\vsp 19& & & & & &  &  \it 12 & \it 6063\\
\cline{1-9}
\vsp 20& & & & & &  &  \it 3 & \it 3540\\
\cline{1-9}
\vsp 21& & & & & &  &  \it 3 & \it 1284\\
\cline{1-9}
\vsp 22& & & & & &  & & \it 761  \\
\cline{1-9}
\vsp 23& & & & & &  & &  \it 124 \\
\cline{1-9}
\vsp 24& & & & & &  & 1 & \it 132\\
\cline{1-9}
\vsp 25& & & & & &  & & \it 39  \\
\cline{1-9}
\vsp 26& & & & & &  & & \it 3  \\
\cline{1-9}
\vsp 27& & & & & &  & &   \\
\cline{1-9}
\vsp 28& & & & & &  & &  1  \\
\cline{1-9}
\vsp Subtotal& 1& 1& 3& 8& {29}& {240}& {2335}&\it 50327\\
\cline{1-9}
\end{tabular}\\
}
\end{table}

\clearpage
\section{Minimal grid diagrams of prime knots with arc index 12 up to 16 crossings}
In the following lists, we used the DT name~\cite{knotinfo}.

In the grid diagrams, we used the convention that the vertical line segment  crosses over the horizontal line segment at every crossing.

\setlength{\unitlength}{0.4pt}

\bigskip
\subsection*{Prime knots with arc index 12 and crossing number 10}~\newline

\bigskip
\noindent
%
 \

\bigskip
\noindent\subsection*{Prime knots with arc index 12 and crossing number 12} ~\newline

\medskip\noindent
%
 \

\clearpage

\noindent\input index12_cr13b

\bigskip
\noindent\input index12_cr14b

\clearpage
\noindent\input index12_cr15b

\noindent\input index12_cr15c

\bigskip
\noindent\input index12_cr16b

\noindent\input index12_cr16c

\noindent\input index12_cr16d

\clearpage

\vfill
\subsection*{Acknowledgments} ~\newline
This work was supported in part by the National Research Foundation of Korea
Grant funded by the Korean Government (NRF-2013-056086).
This work was accelerated by the parallel computing service of KISTI Supercomputing Center.

\end{document}